\documentclass[12pt,a4paper]{amsart}
\usepackage{a4}
\usepackage[english]{babel}
\usepackage{amsmath,amssymb,amsthm}
\def\leq {\leqslant}
\def\le {\leqslant}
\def\ge {\geqslant}
\def\geq {\geqslant}
\def\@bibitem[#1]#2{\item\@biblabel{#1}.\if@filesw
{\def\protect##1{\string##1\space}\immediate\write
\@auxout{\string\bibcite{#2}{#1}}}\fi\ignorespaces\@showtag{#2}}
\newcommand{\no}{\symbol{252}\;}
\theoremstyle{plain}
\newtheorem{theorem}{Theorem}%[section]
\newtheorem{lemma}{Lemma}%[section]
\renewcommand{\theequation}%
{\arabic{section}.\arabic{equation}}
\begin{document}

\title {On orders of approximation functions of generalized smoothnes in  Lorentz spaces}
%\footnote{ }}

\author{ G. Akishev}

\address{ Department of Mathematics and Information Technology,
Buketov Karaganda State University,
Universytetskaya~28 ,
100028, Karaganda , Republic Kazakhstan
\\
Ural Federal University , pr. Lenina 51, \\ Yekaterinburg , 620000, Russia}

\maketitle

\begin{quote}
\noindent{\bf Abstract.}
This paper considers the Lorentz space with mixed norm of periodic functions of many variables and of the generalized Nikol'skii -- Besov classes. Estimates for the order of approximation of the generalized Nikol'skii -- Besov classes by partial sums of Fourier's series for multiple trigonometric system in Lorentz spaces with mixed norm are obtained.
\end{quote}

\vspace*{0.2 cm}



\section{Introduction}
\label{sec1}

Let  $\overline{x} =\left( x_{1},...,x_{m} \right) \in {\Bbb{I}}^{m}
=\left[ 0,2\pi \right]^{m} $  and let $\theta _{j}, p_{j} \in
\left[ 1,+\infty \right)$, $j=1,...,m$, $\Bbb{N}$ be the set of natural numbers.

We shall denote by $ L_{\overline{p},\overline{\theta}}(\Bbb{I}^m)$
the Lorentz spaces with mixed norm  of Lebesgue measurable
functions $f(\overline{x})$ defined on $\Bbb{R}^{m}$ with of period $2\pi$ for each variable such
that
$$
\|f\|_{\overline{p},\overline{\theta}}=\|...\|f\|_{p_{1},\theta_{1}}...\|_{p_{m},\theta_{m}}
< +\infty,
$$
where
$$
\|g\|_{p,\theta}=
\left\{\int\limits_{0}^{2\pi}(g^{*}(t))^{\theta}t^{\frac{\theta}{p}-1}
dt\right\}^{\frac{1}{\theta}},
$$
where $g^{*}$ is a non-increasing rearrangement of the function
$|g|$ (see \cite{ref12}).

As we know, that in case when $p_{j}=\theta_{j}$, $j=1,...,m$, the space $ L_{\overline{p},\overline{\theta}}(\Bbb{I}^m)$
coincides with the Lebesgue space $ L_{\overline{p}}(I^m)$ with mixed norm (for the definition
see \cite{ref21}, p. 128):
 $$
\|f\|_{\bar{p}}=\Biggl[\int_{0}^{2\pi}
\biggl[\cdots\biggl[\int_{0}^{2\pi}|f(\bar x)|
^{p_{1}}dx_{1}\biggr]^{\frac{p_{2}}{p_{1}}}
\cdots\biggr]^{\frac{p_{m}}{p_{m-1}}}dx_{m}\Biggr]^{\frac{1}{p_{m}}}.
$$

Let $\overset{\circ \;\;}{L_{\overline{q} ,\overline{\theta } }}
\left(\Bbb{I}^{m} \right)$ be the set of functions $f\in
L_{\overline{q} ,\overline{\theta}}\left(\Bbb{I}^{m}\right)$ such
that
$$
\int\limits_{0}^{2\pi }f\left(\overline{x}\right) dx_{j}
=0,\;\;\forall j=1,...,m,
$$
and let $a_{\overline{n} } (f)$ be the Fourier coefficients of the function $f\in
L_{1} \left(\Bbb{I}^{m}\right)$ with respect to the multiple
trigonometric system $\{e^{i\langle\overline{n}, \overline{x}\rangle}\}_{\bar{n} \in \Bbb{Z}^{m}}$. Then, we set
$$
\delta _{\overline{s} } \left( f,\overline{x} \right)
=\sum\limits_{\overline{n} \in \rho \left( \overline{s} \right)
}a_{\overline{n} } \left( f\right) e^{i\langle\overline{n}
,\overline{x}\rangle } ,
$$
where $\langle\bar{y},\bar{x}\rangle=\sum\limits_{j=1}^{m}y_{j}
x_{j}$,
 $\rho (\bar{s})=\left\{ \overline{k} =\left( k_{1}
,...,k_{m} \right) \in \Bbb{Z}^{m}: \quad 2^{s_{j} -1} \leq \left| k_{j}
\right| <2^{s_{j} } ,j=1,...,m\right\}$.

A function  $\Omega (\bar t) = \Omega (t_{1},...,t_{m})$ is a function of mixed module continuity type of an order $l\in \Bbb{N}$ if it satisfies the following conditions:

 1) $\Omega (\bar t) > 0$, $t_{j} > 0$, $j=1,...,m$, $\Omega (\bar t) = 0$, if $\prod_{j=1}^{m}t_{j}=0$;

 2) $\Omega (\bar t)$ increases in each variable;

 3) $\Omega (k_{1}t_{1},...,k_{m}t_{m}) \leq \left (\prod_{j=1}^{m}
k_{j}\right)^{l}\Omega (t_{1},...,t_{m})$, $k_{j}\in \Bbb{N}$, $j=1,...,m$;

 4) $\Omega (\bar t)$ is continuous for $t_{j}>0$, $j=1,...,m$.

Let us consider the following sets
 $$
\Gamma (\Omega , N) = \left\{\bar s =(s_{1},...,s_{m})\in
\Bbb{Z}_{+}^{m} : \Omega (2^{-s_{1}},...,2^{-s_{m}})\geq \frac{1}{N}
\right\},
 $$
 $$
Q(\Omega , N)= \cup_{{}_{\bar{s}\in \Gamma (\Omega ,
N)}}\rho(\bar{s}),
 $$
 $$
\Gamma^{\perp} (\Omega , N) =
\Bbb{Z}_{+}^{m}\setminus \Gamma (\Omega ,N), \eqno (1)
 $$
 $$
\Lambda (\Omega , N) = \Gamma^{\perp} (\Omega , N)\setminus \Gamma^{\perp} (\Omega , 2^{l} N). \eqno (2)
 $$

It follows from (1) and (2) that $\Lambda (\Omega , N) \subset \Gamma^{\perp}(\Omega , N)$ and
 $$
\frac{1}{2^{l}N} \le \Omega (2^{- \bar{s}}) < \frac{1}{N} \eqno (3)
 $$
for $\bar{s}\in \Lambda (\Omega , N).$
In  \cite{ref23}, N.N. Pustovoitov proved that $\Lambda(\Omega , N) \neq \emptyset$  and
$$
|\Lambda(\Omega , N)| \asymp \left(\log_{2} N \right)^{m-1}, \eqno(4)
$$
where $|F|$ is the number of elements of the set $F.$

We will use the notation $S_{Q(\Omega , N)}(f,\bar{x})=\sum_{\bar{k}\in Q(\Omega , N)}a_{\bar k}(f)\cdot
e^{i\langle\bar{k},\bar{x}\rangle}$ for a partial sum of the Fourier series of a function $f$.

For a sequence of numbers we write $\left\{ a_{\overline{n} } \right\}
_{\overline{n} \in \Bbb{Z}^{m} } \in l_{\overline{p} }$ if
$$ \left\|
\left\{ a_{\overline{n} } \right\} _{\overline{n} \in \Bbb{Z}^{m} }
\right\| _{l_{\overline{p} } } =\left\{ \sum\limits_{n_{m}
=-\infty }^{\infty }\left[ ...\left[ \sum\limits_{n_{1} =-\infty
}^{\infty }\left| a_{\overline{n} } \right| ^{p_{1} }  \right]
^{\frac{p_{2} }{p_{1} } } ...\right] ^{\frac{p_{m} }{p_{m-1} } }
\right\} ^{\frac{1}{p_{m} } } <+\infty ,
$$
where $\overline{p} =\left( p_{1} ,...,p_{m} \right)$, $1\leq
p_{j} <+\infty$, $j=1,2,...,m$.

For a given function of mixed module smoothness type
 $\Omega(\bar t)$ consider the generalized Nikol'skii -- Besov class
$$
S_{\overline{p} ,\overline{\theta } ,\overline{\tau }
}^{\Omega}B =\left\{ f\in \overset{\circ
\;\;}{L_{\overline{p} ,\overline{\theta } }} \left(I^{m}\right):
\left\| \left\{
\Omega^{-1} (2^{-\bar s}) \left\| \delta
_{\overline{s} } \left( f\right) \right\| _{\overline{p}
,\overline{\theta } } \right\}_{\bar{s} \in \Bbb{Z}_{+}^{m} } \right\| _{l_{\overline{\tau } } }
\leq 1\right\},
$$
where $\overline{p} =\left(p_{1},...,p_{m} \right)$,
$\overline{\theta } = \left( \theta _{1} ,...,\theta _{m}
\right)$, $\overline{\tau } =\left( \tau _{1} ,...,\tau _{m}
\right)$, $1 < p_{j} < +\infty, $ $1<\theta _{j} < \infty,$ $1 \le \tau_{j} \le +\infty$,
$j=1,...,m$, and $\Omega (2^{-\bar s}) = \Omega (2^{-s_{1}},...,2^{-s_{m}}).$

If  $\Omega(\bar t) = \prod_{j=1}^{m}t_{j}^{r_{j}},$ $\,\, r_{j}>0, \,\, j=1,...,m$, then this class is denoted by $S_{\overline{p} ,\overline{\theta },\overline{\tau}}^{\overline{r}}B.$

In case $ p_{j} = \theta_{j}=p$ and $\Omega(\bar t) = \prod_{j=1}^{m}t_{j}^{r_{j}}, r_{j} < l, \tau_{j}=+\infty, j=1,...,m,$
$S_{\overline{p} ,\overline{\theta} ,\overline{\tau}}^{\Omega}B$ was defined by S.M. Nikol'skii \cite{ref19}, and for  $1\le \tau_{j} < +\infty, j=1,...,m$, by T.I. Amanov \cite{ref06} and P.I. Lizorkin, S.M. Nikol'skii \cite{ref18}.

As pointed out in \cite{ref35}, \cite{ref36} one of the difficulties in the theory of
approximation of functions of several variables is the choice of harmonics of the approximating polynomials. The first author, who
suggested to approximate functions of several variables by
polynomials with harmonics in hyperbolic crosses, was K.I. Babenko \cite{ref07}. After that approximations of various classes of smooth
functions by this method were considered by S.A. Telyakovskii \cite{ref32}, B.S. Mityagin  \cite{ref19}, Ya. S. Bugrov \cite{ref13}, N.S. Nikol'skaya \cite{ref22}, E.M. Galeev \cite{ref16, ref17}, V.N. Temlyakov \cite{ref33, ref34}, Dinh Dung \cite{ref15}, A.R. DeVore, S.V. Konyagin and V.N. Temlyakov \cite{ref14}, H. - J. Schmeisser and W. Sickel \cite{ref29}, W. Sickel and T. Ullrich \cite{ref27},
  A.S. Romanyuk \cite{ref25, ref26}.

For the generalized Besov class this problem was considered by N.N. Pustovoitov \cite{ref23}, \cite{ref24},
Sun Yongsheng and Wang Heping \cite{ref31}, D.B. Bazakhanov \cite{ref09}, M. Sikhov \cite{ref28}, and S.A. Stasyuk \cite{ref30}.

Exact orders of the approximation of the Nikol'skii--Besov classes in the metric of the Lorentz space were found by the author \cite{ref01, ref02} and K.A. Bekmaganbetov \cite{ref10}, \cite{ref11}.

An order of approximation of the class
$S_{\overline{p} ,\overline{\theta },\overline{\tau}}^{\overline{r}}B$ by partial Fourier sums
 $S_{n}^{\bar \gamma}(f,\bar{x})=\sum_{\langle\bar{s}, \bar\gamma \rangle < n }\delta_{\overline{s}} \left( f, \bar x\right)$
was found in \cite{ref01}. In \cite{04} for class $S_{\overline{p} ,\overline{\theta } ,\overline{\tau}
}^{\Omega}B$ proved following statement.

{\bf Theorem} (see \cite{ref04}).
Let  $1\le \theta_{j}^{(1)}, \theta_{j}^{(2)}, \tau_{j} < +\infty,$ $1 < p_{j} < q_{j} < \infty,$ $j=1,...,m,$ and $\Omega(\bar t)$ be a function of mixed module continuity type of an order  $l$, which satisfies the conditions
 $(S)$ and $(S_{l})$, $\alpha_{j} > \frac{1}{p_{j}} - \frac{1}{q_{j}},$ $j=1,...,m.$

1) If  $1\le \theta_{j}^{(2)}< \tau_{j} < +\infty,$ $j=1,...,m,$ then
$$
\frac{1}{N}\left(\log_{2}N\right)^{-\sum\limits_{j=2}^{m}
\frac{1}{\tau_{j}}}\left\|\left\{\prod_{j=1}^{m}2^{s_{j}\left( \frac{1}{p_{j}} - \frac{1}{q_{j}}\right)}  \right\}_{\bar s \in \Lambda(N)}\right\|_{l_{\bar\theta^{(2)}}} <<
\sup \limits_{f\in
S_{\overline{p},\bar{\theta}^{(1)} \overline{\tau}
}^{\Omega}B }\|f-S_{Q(N)}(f)\|_{\overline{q}, \overline{\theta}^{(2)}} <<
$$
$$
<< \frac{1}{N} \left\|\left\{\prod_{j=1}^{m}2^{s_{j}\left( \frac{1}{p_{j}} - \frac{1}{q_{j}}\right)}  \right\}_{\bar s \in \Lambda(N)}\right\|_{l_{\bar\epsilon}},
$$
where $\bar\epsilon = (\epsilon_{1},...,\epsilon_{m}),$ $\epsilon_{j} = \frac{\tau_{j}\theta_{j}^{(2)}}{\tau_{j}-\theta_{j}^{(2)}}, j=1,...,m.$

2) If $\tau_{j}\le \theta_{j}^{(2)} , j=1,...,m,$ then
$$
\sup \limits_{\bar s \in \Lambda(N)}
\Omega(2^{-\bar s})
\prod_{j=1}^{m}2^{s_{j}\left( \frac{1}{p_{j}} - \frac{1}{q_{j}}\right)} <<
\sup \limits_{f\in
S_{\overline{p},\bar{\theta}^{(1)} \overline{\tau}
}^{\Omega}B }\|f-S_{Q(N)}(f)\|_{\overline{q}
,\overline{\theta}^{(2)}} <<
$$
$$
<< \sup \limits_{\bar s \in \Gamma^{\perp}(N)}
\Omega(2^{-\bar s})
\prod_{j=1}^{m}2^{s_{j}\left( \frac{1}{p_{j}} - \frac{1}{q_{j}}\right)}.
$$

The notation $A\left( y\right) \asymp
B\left( y\right)$ means that there exist positive constants
 $C_{1},\,C_{2} $ such that  $C_{1}A\left( y\right)
\leq B\left( y\right) \leq C_{2}A\left( y\right)$. If $B \le C_{2}A$ or $A\ge C_{1}B$, then we write $B<<A$ or $A>>B$.

The main aim of the present paper is to estimate the order of
the quantity
$$
\sup \limits_{f\in
S_{\overline{p} ,\overline{\theta } ,\overline{\tau }
}^{\Omega}B }\|f-S_{Q(\Omega, N)}(f)\|_{\overline{q}
,\overline{\theta}} .
$$

 This paper is organized as follows. In the
second section some auxiliary lemmas are given.
The third section establishes the estimate of the order approximation of the
Nikol'skii--Besov classes in the Lorentz space with mixed norm.

\vspace*{0.2 cm}

\section {Auxiliary lemmas}

\vspace*{0.2 cm}

In what follows, we denote by
 $\chi_{\varkappa(n)}(\bar s)$ the characteristic function of the set
 $ \varkappa(n)=\{\bar{s}=(s_{1},...,s_{m})\in\mathbb{Z}_{+}^{m}:
\quad \langle\bar{s},\bar{\gamma}\rangle=n\} $.

\begin{lemma} %Лемма 1
{\it Let
 $\bar\tau=(\tau_{1},...,\tau_{m}) $, $1\leq \tau_{j}< +\infty$, $j=1,...,m$. Then the following relation holds:
$$
\left\|\left\{\chi_{\varkappa(n)}(\bar s)\right\}_{\bar{s}
\in\varkappa(n)}\right\|_{l_{\bar\tau}}\asymp
n^{\sum\limits_{j=2}^{m}\frac{1}{\tau_{j}}}.
$$}
\end{lemma}
\begin{lemma} %Лемма 2
{\it Let $\bar\gamma=(\gamma_{1},...,\gamma_{m})$, $\overline{\gamma}^{'} =\left( \gamma_{1}^{'}
,...,\gamma_{m}^{'}\right)$, $\gamma_{j}^{'}=\gamma _{j},
j=1,...,\nu $, $1<\gamma_{j}<\gamma_{j}^{'} , j=\nu+1,...,m$, and
let $\bar\tau=(\tau_{l+1},...,\tau_{m}), $ where $1\leq \tau_{j}<
+\infty$, $j=1,...,m$,
and $\alpha >0$. Then the following relation holds:
$$
 I_{n}^{(l)} = \left\|\left\{2^{-\alpha \langle\bar s,\bar \gamma^{'} \rangle} \right\}_{\bar{s}
\in\varkappa(n)}\right\|_{l_{\bar\tau}}\asymp
2^{-n\alpha}\cdot n^{\sum\limits_{j=2}^{\nu}\frac{1}{\tau_{j}}}.
$$}
\end{lemma}
Lemma 1, 2 are proved in \cite{ref02}.

Let us recall definitions of the conditions $(S),$ $(S_{l})$ given by S.B.Stechkin and N.K. Bary  \cite{ref08}.

 {\bf Definition.} A function $g(t)$ satisfies the condition $(S)$,
if for some  $\alpha \in (0 , 1)$ the function  $t^{-\alpha}g(t)$
almost increases on $(0 , 1].$

We say that a function  $\Omega (\bar t)$ satisfies the condition $(S)$ on  $(0 ,1]^{m},$ if it satisfies this condition on each variable.

{\bf Definition.} A function $g(t)$ satisfies the condition  $(S_{l}),$
if for some   $\alpha \in (0 , l)$ the function $t^{-\alpha}g(t)$
almost decreases on  $(0 , 1].$

We say that a function   $\Omega (\bar t)$ satisfies the condition $(S_{l})$ on $(0 ,1]^{m},$ if it satisfies this condition on each variable.

\begin{lemma} (see \cite{ref04}) %Лемма 3
 Let  $1\leq \theta _{j} < +\infty$, $j=1,...,m$, and $\Omega (\bar
t)$ be a function of mixed module continuity type of an order $l$
which satisfies the $(S)$-condition for  $\bar \alpha =
(\alpha_{1},...,\alpha_{m})$, $\alpha_{j} > \beta_{j} \ge 0$, $j =
1,...,m$. Then for  $1 \le \theta_{j} < +\infty$, $j=1,...,m$, the following relation holds
$$
\left\| \left\{\Omega (2^{-s_{1}},...,2^{-s_{m}})
\prod_{j=1}^{m}2^{s_{j}\beta_{j}} \right\}_{\overline{s} \in
\Gamma^{\perp} (\Omega , N)}
\right\|_{l_{\overline{\theta}}}\asymp
$$
$$
\asymp \left\| \left\{\Omega (2^{-s_{1}},...,2^{-s_{m}})
\prod_{j=1}^{m}2^{s_{j}\beta_{j}} \right\} _{\bar{s} \in
\Lambda(N)} \right\|_{l_{\bar{\theta}}}.
$$
\end{lemma}

\begin{lemma} (see \cite{ref04}). %Лемма4
Let $\Omega(\bar t)$ be a function of mixed module continuity type of an order  $l$, which satisfies the conditions
 $(S)$ and $(S_{l})$, $1\le \tau_{j} < +\infty, j=1,...,m$, and $\Lambda(\Omega, N) =\Gamma^{\perp}(\Omega, N)\setminus\Gamma^{\perp}(\Omega, 2^{l}N).$ Then
$$
\left\|\left\{\chi_{\Lambda(\Omega, N)}(\bar s) \right\}_{\bar s \in \Lambda(\Omega, N)}\right\|_{l_{\bar{\tau}}}\asymp \left(\log_{2} N \right)^{\sum\limits_{j=2}^{m}\frac{1}{\tau_{j}}}.
$$
\end{lemma}

{\bf Remark.} Note that for the case  $\tau_{1}=...=\tau_{m}=1$ Lemma 4 was proved by N.N. Pustovoitov \cite{ref23}.

\begin{theorem} %Теорема 1
Let $\bar{q}=(q_{1},...,q_{m})$, $1<q_{j}<\infty$, $ j=1,...,m$,
$\beta = min \{q_{1},...,q_{m},2 \}.$ Then, for any function  $f
\in L_{\bar q}(I^{m})$, the following inequality holds
$$
\|f\|_{\bar q} << \left\{\sum\limits_{\bar s \in
\mathbb{Z}_{+}^{m}} \| \delta_{\bar s} (f)\|_{\bar
q}^{\beta}\right\}^{\frac{1}{\beta}}.
$$
\end{theorem}
The proof of theorem is given in \cite{ref03}.

\begin{theorem} (see \cite{ref01}).  %Теорема 2
{\it Let
$ \bar p =(p_{1},...,p_{m})$, $\bar q =(q_{1},...,q_{m})$,
$\bar\theta^{(1)} =(\theta_{1}^{(1)},...,\theta_{m}^{(1)})$,
$\bar\theta^{(2)} =(\theta_{1}^{(2)},...,\theta_{m}^{(2)})$.
Assume that  $1\leq p_{j}<q_{j}<+\infty$, $1\leq\theta_{j}^{(1)}, \theta_{j}^{(2)}  <+\infty$, $j=1,...,m$.
If $f\in \overset{\circ \; \;}{L}_{\bar
p,\bar\theta^{(1)}}(\Bbb{I}^{m})$,
$\max_{j=1,...,m-1}{\theta_{j}^{(2)}} < \min_{j=2,...,m}{q_{j}}$
and the quantity
$$
\sigma(f)\equiv\Biggl\{\sum_{s_{m}=1}^{\infty}
 2^{s_{m}\theta_{m}^{(2)}(\frac{1}{p_{m}}-\frac{1}{q_{m}})}\biggl[\cdots\biggl[
\sum_{s_{1}=1}^{\infty}2^{s_{1}\theta_{1}^{(2)}(\frac{1}{p_{1}}-\frac{1}{q_{1}})}
\|\delta_{\bar s}(f)\|_{\bar
p,\bar\theta^{(1)}}^{\theta_{1}^{(2)}}\biggr]^{
\frac{\theta_{2}^{(2)}}{\theta_{1}^{(2)}}}
\cdots\biggr]^{\frac{\theta_{m}^{(2)}}{\theta_{m-1}^{(2)}}}
\Biggr\}^{\frac{1}{\theta_{m}^{(2)}}}
$$
is finite, then $f\in
\overset{\circ\;\;}{L_{\bar{q},\bar\theta^{(2)}}}(\Bbb{I}^{m})$ and
$$
\|f\|_{\bar{q},\bar\theta^{(2)}} << \sigma(f).
$$}
\end{theorem}

\begin{theorem} (see \cite{ref01}). %Теорема 3
{\it Let $\bar q =(q_{1},...,q_{m})$ ,
$\bar\theta =(\theta_{1},...,\theta_{m})$, $\bar\lambda = (\lambda_{1},...,\lambda_{m})$.
 Assume that
 $1<q_{j}<\tau_{j}<+\infty$, $1<\theta_{j}<+\infty$, $j=1,...,m$.
If $f\in \overset{\circ\;\;}{L}_{\bar{q},\bar{\theta}}(\Bbb{I}^{m})$ and
$$
f(\bar x)\sim \sum_{\bar{s}\in\mathbb{Z}_{+}^{m}}b_{\bar s}
\sum_{\bar{k}\in\rho(\bar s)}e^{i\langle\bar{k},\bar{x}\rangle} ,
$$
then
$$
\|f\|_{\bar{q},\bar{\theta}} >> \left\{\sum_{s_{m}=1}^{\infty}2^{s_{m}\theta_{m}\left(\frac{1}{\lambda_{m}}-
\frac{1}{q_{m}}\right)}\left[...\left[
\sum_{s_{1}=1}^{\infty}2^{s_{1}\theta_{1}\left(\frac{1}{\lambda_{1}}-
\frac{1}{q_{1}}\right)}\left(\|\delta_{\bar
s}(f)\|_{\bar{\lambda},\bar{\theta}}
\right)^{\theta_{1}}\right]^{\frac{\theta_{2}}{\theta_{1}}}...
\right]^{\frac{\theta_{m}}{\theta_{m-1}}}\right\}^{\frac{1}{\theta_{m}}}.
$$}
\end{theorem}

\vspace*{0.2 cm}
\begin{center}
\section {Main results}
\end{center}
\vspace*{0.2 cm}

Let us prove the main results of the present paper.

 Consider the function $\Omega_{1}(\bar{t}) = \Omega(\bar{t})\prod\limits_{j=1}^{m}t_{j}^{-(\frac{1}{p_{j}} - \frac{1}{q_{j}})}$ для $t_{j} \in (0, 1], j = 1,...,m$ and respectively set $Q(\Omega_{1}, N), \Gamma^{\perp}(\Omega_{1}, N), \Lambda(\Omega_{1}, N)$.

\begin{theorem} %теорема 4
Let  $1\le \theta_{j}^{(1)}, \theta_{j}^{(2)}, \tau_{j} < +\infty,$ $1 < p_{j} < q_{j} < \infty,$ $j=1,...,m,$ and $\Omega(\bar t)$ be a function of mixed module continuity type of an order  $l$, which satisfies the conditions
 $(S)$ and $(S_{l})$, $\alpha_{j} > \frac{1}{p_{j}} - \frac{1}{q_{j}},$ $j=1,...,m$ и $\Omega_{1}(\bar{t}) = \Omega(\bar{t})\prod\limits_{j=1}^{m}t_{j}^{-(\frac{1}{p_{j}} - \frac{1}{q_{j}})}$.

1) If  $1\le \theta_{j}^{(2)}< \tau_{j} < +\infty,$ $j=1,...,m,$ then
$$
\sup \limits_{f\in
S_{\overline{p},\bar{\theta}^{(1)} \overline{\tau}
}^{\Omega}B }\|f-S_{Q(\Omega_{1}, N)}(f)\|_{\overline{q}, \overline{\theta}^{(2)}} \asymp \frac{1}{N}(\log_{2}N)^{\sum\limits_{j=2}^{m}(\frac{1}{\theta_{j}^{(2)}} - \frac{1}{\tau_{j}})}.
$$

2) If $\tau_{j}\le \theta_{j}^{(2)} , j=1,...,m,$ then
$$
\sup \limits_{f\in
S_{\overline{p},\bar{\theta}^{(1)} \overline{\tau}
}^{\Omega}B }\|f-S_{Q(\Omega_{1}, N)}(f)\|_{\overline{q}
,\overline{\theta}^{(2)}} \asymp \frac{1}{N}.
$$
\end{theorem}

{\bf Proof.} Taking into account $\delta_{\bar{s}}(f - S_{Q(\Omega_{1}, N)}(f)) = 0$
, if $\bar{s} \in Q(\Omega_{1}, N)$ and $\delta_{\bar{s}}(f - S_{Q(\Omega_{1}, N)}(f)) = \delta_{\bar{s}}(f)$
, if $\bar{s} \notin Q(\Omega_{1}, N)$ By Theorem 2, we have
$$
\|f-S_{Q(\Omega_{1}, N)}(f)\|_{\overline{q}
,\overline{\theta}^{(2)}} <<
\left\|\left\{\prod_{j=1}^{m}2^{s_{j}\left( \frac{1}{p_{j}} - \frac{1}{q_{j}}\right)}\left\|\delta_{\bar{s}}(f-S_{Q(\Omega_{1}, N)}(f))  \right\|_{\overline{p}, \bar{\theta}^{(1)}}\right\}_{\bar{s}\in \Bbb{Z}_{+}^{m}}\right\|_{l_{\bar{\theta}^{(2)}}} =
$$
$$
= C\left\|\left\{\prod_{j=1}^{m}2^{s_{j}\left( \frac{1}{p_{j}} - \frac{1}{q_{j}}\right)}\left\|\delta_{\bar{s}}(f-S_{Q(\Omega_{1}, N)}(f))  \right\|_{\overline{p}, \bar{\theta}^{(1)}}\right\}_{\bar{s}\in \Gamma^{\perp}(\Omega_{1}, N)}\right\|_{l_{\bar{\theta}^{(2)}}}.
$$
for any function  $f\in S_{\overline{p},\bar{\theta}^{(1)}, \overline{\tau}}^{\Omega}B$.

Since  $\beta_{j} = \frac{\tau_{j}}{\theta_{j}^{(2)}} > 1, j=1,...,m,$ and by applying Holder's inequality we obtain the following
$$
\|f-S_{Q(\Omega_{1}, N)}(f)\|_{\overline{q}
,\overline{\theta}^{(2)}} <<
\left\|\left\{\Omega^{-1}(2^{-\bar{s}})\left\|\delta_{\bar{s}}(f)  \right\|_{\overline{p}, \bar{\theta}^{(1)}}\right\}_{\bar{s}\in \Bbb{Z}_{+}^{m}}\right\|_{l_{\bar{\tau}}}\times
$$
$$
\times\left\|\left\{\Omega(2^{-\bar{s}})\prod_{j=1}^{m}2^{s_{j}\left( \frac{1}{p_{j}} - \frac{1}{q_{j}}\right)}\right\}_{\bar{s}\in \Gamma^{\perp}(\Omega_{1}, N)}\right\|_{l_{\bar{\epsilon}}},   \eqno (5)   %(13) %(*1)
$$
where $\bar\epsilon = (\epsilon_{1},...,\epsilon_{m}),$ $\epsilon_{j} = \frac{\tau_{j}\theta_{j}^{(2)}}{\tau_{j}-\theta_{j}^{(2)}},$ $j=1,...,m.$

2.	Since by assumption theorem  the function $\Omega(\bar{t})$ satisfies $S$ and  $S_{l}$ conditions and $\alpha_{j} > \frac{1}{p_{j}} - \frac{1}{q_{j}},$ $j=1,...,m$ , then the function $\Omega_{1}(\bar{t})$ satisfies conditions $S$ and  $S_{l}$.

Therefore , by Lemma 3, Lemma 4 and the definition of the set $\Gamma^{\perp}(\Omega_{1}, N)$ in (5), we have
$$
\sup \limits_{f\in
S_{\overline{p},\bar{\theta}^{(1)}, \overline{\tau}
}^{\Omega}B }\|f-S_{Q(\Omega_{1}, N)}(f)\|_{\overline{q}
,\overline{\theta}^{(2)}} <<
\left\|\left\{\Omega_{1}(2^{\bar{s}})\right\}_{\bar{s}\in \Lambda(\Omega_{1}, N)}\right\|_{l_{\bar{\epsilon}}} <<
$$
$$
<< \frac{1}{N}\left\|\left\{\chi_{\Lambda(\Omega_{1}, N)}(\bar s) \}_{\bar{s}\in \Lambda(\Omega_{1}, N)}\right\}\right\|_{l_{\bar{\epsilon}}}<<
\frac{1}{N}(\log_{2}N)^{\sum\limits_{j=2}^{m}(\frac{1}{\theta_{j}^{(2)}} - \frac{1}{\tau_{j}})}.
$$

In item 1) of the theorem the upper bound has been proved.

Let us prove the lower bound.
Consider the function
$$
f_{0}(\bar x)=\left(\log_{2} N\right)^{-\sum\limits_{j=2}^{m}\frac{1}{\tau_{j}}}
\sum_{\bar{s}\in \Lambda(\Omega_{1}, N)}\prod_{j=1}^{m}\Omega(2^{-\bar s})
2^{-s_{j}\left(1-\frac{1}{p_{j}}\right)}\sum_{\bar{k}\in\rho(\bar
s)} e^{i\langle\bar{k},\bar{x}\rangle}.
$$

In one-dimensional case for Dirichlet's kernel
$D_{n}(x) = \frac{1}{2} + \sum\limits_{k=1}^{n}e^{ikx}$
the following statement holds
 $$
\|D_{n}\|_{p, \theta}
\asymp
n^{1-\frac{1}{p}}, \;\;  1< p < +\infty, \;\; 1< \theta < +\infty.
 $$

Then, by the property of the norm, we have
$$
\left\|\sum_{k_{j} = 2^{s_{j} - 1}}^{2^{s_{j}} - 1} e^{ik_{j}x_{j}}\right\|_{p_{j}, \theta_{j}^{(1)}} \le \left\|D_{2^{s_{j}} - 1} \right\|_{p_{j}, \theta_{j}^{(1)}} + \left\|D_{2^{s_{j} - 1} - 1} \right\|_{p_{j}, \theta_{j}^{(1)}} << 2^{s_{j}(1 - \frac{1}{p_{j}})},
$$
provided $1 < p_{j} < +\infty$, $1 < \theta_{j}^{(1)} < +\infty$, $j =1,...,m$.
Hence
$$
\left\|\sum_{\bar{k}\in\rho(\bar
s)} e^{i\langle\bar{k},\bar{x}\rangle}\right\|_{\overline{p}, \bar{\theta}^{(1)}} = \prod_{j = 1}^{m} \left\|\sum_{k_{j} = 2^{s_{j} - 1}}^{2^{s_{j}} - 1} e^{ik_{j}x_{j}}\right\|_{p_{j}, \theta_{j}^{(1)}} << \prod_{j=1}^{m}2^{s_{j}\left(1-\frac{1}{p_{j}}\right)}.
$$
Let us prove the rest of the equality. By Lemma B in [1], the following inequality holds
$$
\max_{\bar{x} \in I^{m}} |\sum_{\bar{k}\in\rho(\bar
s)} e^{i\langle\bar{k},\bar{x}\rangle}| << \prod_{j=1}^{m}2^{s_{j}\left(1-\frac{1}{p_{j}}\right)}
\left\|\sum_{\bar{k}\in\rho(\bar
s)} e^{i\langle\bar{k},\bar{x}\rangle}\right\|_{\overline{p}, \bar{\theta}^{(1)}}.       \eqno (6) %(14)
$$
It is known that
$$
\max_{\bar{x} \in I^{m}} |\sum_{\bar{k}\in\rho(\bar
s)} e^{i\langle\bar{k},\bar{x}\rangle}| \ge |\sum_{\bar{k}\in\rho(\bar
s)} e^{i\langle\bar{k},\bar{0}\rangle}| \ge 2^{-m}\prod_{j=1}^{m}2^{s_{j}}.
$$
Therefore, it follows from (6) that
$$
 \prod_{j=1}^{m}2^{s_{j}\left(1-\frac{1}{p_{j}}\right)} << \left\|\sum_{\bar{k}\in\rho(\bar
s)} e^{i\langle\bar{k},\bar{x}\rangle}\right\|_{\overline{p}, \bar{\theta}^{(1)}}.
$$
Thus, we have proved the relation
$$
\left\|\sum_{\bar{k}\in\rho(\bar
s)} e^{i\langle\bar{k},\bar{x}\rangle}\right\|_{\overline{p}, \bar{\theta}^{(1)}} \asymp  \prod_{j=1}^{m}2^{s_{j}\left(1-\frac{1}{p_{j}}\right)}.
\eqno (7)
$$
Therefore, by Lemma 4 and by estimation (7), we have
$$
 \left\|\left\{\Omega^{-1}(2^{-\bar{s}})\left\|\delta_{\bar{s}}(f_{0})  \right\|_{\overline{p}, \bar{\theta}^{(1)}}\right\}_{\bar{s}\in \Bbb{Z}_{+}^{m}}\right\|_{l_{\bar{\tau}}} =
$$
$$
= \left\|\left\{\Omega^{-1}(2^{-\bar{s}})\left(\log_{2} N\right)^{-\sum\limits_{j=2}^{m}\frac{1}{\tau_{j}}}
\Omega(2^{-\bar{s}})\prod_{j=1}^{m}2^{- s_{j}\left(1-\frac{1}{p_{j}}\right)}
 \Bigl\|\sum_{\bar{k}\in\rho(\bar
s)} e^{i\langle\bar{k},\bar{x}\rangle}\Bigr\|_{\overline{p}, \bar{\theta}^{(1)}}\right\}_{\bar{s}\in \Lambda(\Omega_{1}, N)}\right\|_{l_{\bar{\tau}}}
$$
$$
= \left(\log_{2} N\right)^{-\sum\limits_{j=2}^{m}\frac{1}{\tau_{j}}}
\left\|\left\{\chi_{\Lambda(\Omega_{1}, N)}(\bar s) \right\}_{\bar s \in \Lambda(\Omega_{1}, N)}\right\|_{l_{\bar\tau}}\le C_{0}.
$$
Hence  $C_{0}^{-1}f_{0} \in S_{\overline{p}, \overline{\tau}}^{\Omega}B.$
Now taking into account that $S_{Q(\Omega_{1}, N)}^{\bar \gamma}(f_{0}, \bar x)= 0, \bar{x} \in I^{m}$ and using Theorem 4 and (7), Lemma 4, we obtain
$$
\|f_{0}-S_{Q(\Omega_{1}, N)}(f_{0})\|_{\overline{q}, \overline{\theta}^{(2)}} = \|f_{0}\|_{\overline{q},\overline{\theta}^{(2)}} >>
$$
$$
 >> \left\|\left\{\prod_{j=1}^{m}2^{s_{j}\left( \frac{1}{\lambda_{j}} - \frac{1}{q_{j}}\right)}\left\|\delta_{\bar{s}}(f_{0})  \right\|_{\overline{\lambda}, \bar{\theta}^{(1)}}\right\}_{\bar{s}\in \Bbb{Z}_{+}^{m}}\right\|_{l_{\bar{\theta}^{(2)}}} =
$$
$$
= C \left\|\left\{\prod_{j=1}^{m}2^{s_{j}\left( \frac{1}{\lambda_{j}} - \frac{1}{q_{j}}\right)}\left(\log_{2} N\right)^{-\sum\limits_{j=2}^{m}\frac{1}{\tau_{j}}} \Omega(2^{-\bar{s}})
\prod_{j=1}^{m}2^{-s_{j}\left(1-\frac{1}{p_{j}}\right)}
 \left\|\sum_{\bar{k}\in\rho(\bar
s)} e^{i\langle\bar{k},\bar{x}\rangle}\right\|_{\overline{\lambda}, \bar{\theta}^{(1)}}\right\}_{\bar{s}\in \Lambda(\Omega_{1}, N) }\right\|_{l_{\bar{\theta}^{(2)}}} >>
$$
$$
>> \left(\log_{2} N\right)^{-\sum\limits_{j=2}^{m}\frac{1}{\tau_{j}}}
 \left\|\left\{\Omega_{1}(2^{-\bar{s}})\right\}_{\bar{s}\in \Lambda(\Omega_{1}, N) }\right\|_{l_{\bar{\theta}^{(2)}}} >>
$$
$$
>> \frac{1}{N}
\left(\log_{2} N\right)^{-\sum\limits_{j=2}^{m}\frac{1}{\tau_{j}}}
 \left\|\left\{\chi_{\Lambda(\Omega_{1}, N)}\right\}_{\bar{s}\in \Lambda(\Omega_{1}, N)}\right\|_{l_{\bar{\theta}^{(2)}}} >>
 \frac{1}{N}
\left(\log_{2} N\right)^{\sum\limits_{j=2}^{m}(\frac{1}{\theta_{j}^{(2)}}-\frac{1}{\tau_{j}})}.
$$
Thus,
$$
\sup \limits_{f\in
S_{\overline{p},\bar{\theta}^{(1)}, \overline{\tau}
}^{\Omega}B }\|f-S_{Q(\Omega_{1}, N)}(f)\|_{\overline{q}
,\overline{\theta}^{(2)}} >> \frac{1}{N}\left(\log_{2} N\right)^{\sum\limits_{j=2}^{m}(\frac{1}{\theta_{j}^{(2)}}-\frac{1}{\tau_{j}})}.
$$
 Item  1) of the theorem has been proved.

Let us prove item 2) of the theorem. Since $\tau_{j}\le \theta_{j}^{(2)}$, $j=1,...,m$, then  by applying Theorem 2 and Jensen's inequality (see \cite{ref21}, p. 125), we obtain
$$
\|f-S_{Q(\Omega_{1}, N)}(f)\|_{\overline{q}, \overline{\theta}^{(2)}}
<< \left\|\left\{\prod_{j=1}^{m}2^{s_{j}\left( \frac{1}{p_{j}} - \frac{1}{q_{j}}\right)}\left\|\delta_{\bar{s}}(f)  \right\|_{\overline{p}, \bar{\theta}^{(1)}}\right\}_{\bar{s}\in \Gamma^{\perp}(\Omega_{1}, N)}\right\|_{l_{\bar{\tau}}} <<
$$
$$
<< \left\|\left\{\Omega^{-1}(2^{-\bar{s}})\left\|\delta_{\bar{s}}(f)  \right\|_{\overline{p}, \bar{\theta}^{(1)}}\right\}_{\bar{s}\in \Bbb{Z}_{+}^{m}}\right\|_{l_{\bar{\tau}}}
\sup\limits_{\bar s \in \Gamma^{\perp}(\Omega_{1}, N)}\Omega(2^{-\bar{s}})\prod_{j=1}^{m}2^{s_{j}\left( \frac{1}{p_{j}} - \frac{1}{q_{j}}\right)} << \frac{1}{N}
$$
for any function $f\in S_{\overline{p},\bar{\theta}^{(1)}, \overline{\tau}}^{\Omega}B$,
which proves the upper bound in item 2). For the lower bound, consider the function
$$
f_{1}(\bar x)= \Omega(2^{-\bar{\widetilde s}})
2^{-\sum\limits_{j=1}^{m}\widetilde{s}_{j}\left(1
-\frac{1}{p_{j}}\right)}\sum_{\bar{k}\in\rho(\bar{\widetilde s})}
e^{i\langle\bar{k},\bar{x}\rangle},
$$
where $\bar{\widetilde{s}} = (\widetilde{s}_{1},...,\widetilde{s}_{m})\in \Lambda(\Omega_{1}, N)$. Then  $f_{1} \in S_{\overline{p},\bar{\theta}^{(1)}, \overline{\tau}}^{\Omega}B.$
Next, by  (7), we have
$$
\|f_{1}-S_{Q(\Omega_{1}, N)}(f_{1})\|_{\overline{q}, \overline{\theta}^{(2)}}=\|f_{1}\|_{\overline{q}, \overline{\theta}^{(2)}} >>
$$
$$
>> \Omega(2^{-\bar{\widetilde s}})
2^{-\sum\limits_{j=1}^{m}\widetilde{s}_{j}\left(1
-\frac{1}{p_{j}}\right)}\prod_{j=1}^{m}2^{\widetilde{s}_{j}\left(1
-\frac{1}{q_{j}}\right)} = C\Omega(2^{-\bar{\widetilde s}})
\prod_{j=1}^{m}2^{\widetilde{s}_{j}\left(\frac{1}{p_{j}}
-\frac{1}{q_{j}}\right)}
$$
for $\bar{\widetilde s} \in \Lambda(\Omega_{1}, N)$.

Hence , by (3) we obtain
$$
\sup \limits_{f\in
S_{\overline{p},\bar{\theta}^{(1)}, \overline{\tau}
}^{\Omega}B }\|f-S_{Q(\Omega_{1}, N)}(f)\|_{\overline{q}
,\overline{\theta}^{(2)}} >> \frac{1}{N}.
$$
This proves the lower bound in item 2).

\begin{theorem} %Теорема 5
Let $\Omega(\bar t)$ be a function of mixed module continuity type of an order  $l$ which satisfies the conditions
 $(S)$ and $(S_{l})$, $1 < q_{j} < p_{j} < \infty$, $p_{j} \ge 2$, $1< \theta_{j}< \infty$, $1\le \tau_{j} \le +\infty,$ $j=1,...,m$.

1) If $2< \tau_{j} < +\infty,$
$j=1,...,m,$  then
$$
\sup \limits_{f\in
S_{\overline{p}, \overline{\tau}
}^{\Omega}B }\|f-S_{Q(\Omega, N)}(f)\|_{\overline{q}
,\overline{\theta}} \asymp \frac{1}{N}\left(\log_{2}N\right)^{\sum\limits_{j=2}^{m}
\left(\frac{1}{2}-\frac{1}{\tau_{j}}\right)}.
$$

2) If $\tau_{j} \le 2$, $j=1,...,m$, then
$$
\sup \limits_{f\in
S_{\overline{p}, \overline{\tau}
}^{\Omega}B }\|f-S_{Q(\Omega, N)}(f)\|_{\overline{q}
,\overline{\theta}} \asymp \frac{1}{N}.
$$

3) If $1 < q_{j} < p_{j} \le 2, \;\; j=1,...,m$ and $p_{0} = \min\{p_{1},...,p_{m}\} < \tau_{j}, \;\; j=1,...,m,$ then
$$
\frac{1}{N}(\log_{2}N)^{\sum\limits_{j=2}^{m}
\left(\frac{1}{p_{j}}-\frac{1}{\tau_{j}}\right)}<< \sup \limits_{f\in
S_{\overline{p}, \overline{\tau}}^{\Omega} B }\|f-S_{Q(\Omega, N)}(f)\|_{\overline{q}, \overline{\theta}} << \frac{1}{N}(\log_{2}N)^{\sum\limits_{j=2}^{m}
\left(\frac{1}{p_{0}}-\frac{1}{\tau_{j}}\right)}.
$$
\end{theorem}

{\bf Proof.}
Since $q_{j} < p_{j}$, $j=1,...,m$, then $L_{\bar{p}}(\Bbb{I}^{m}) \subset L_{\bar{q}, \bar{\theta}}(\Bbb{I}^{m})$ and we have
$$
\|f\|_{\overline{q},\overline{\theta}} << \|f\|_{\overline{p}}, \quad
f\in L_{\bar p}(\Bbb{I}^{m}).
$$
Therefore $S_{\overline{p}, \overline{\tau}}^{\Omega}B \subset  L_{\bar{q}, \bar{\theta}}(\Bbb{I}^{m})$ and
$$
\|f-S_{Q(\Omega, N)}(f)\|_{\overline{q},\overline{\theta}} << \|f-S_{Q(\Omega, N)}(f)\|_{\overline{p}}=
$$
$$
= C\left\|\sum_{\bar{s} \in \Gamma^{\perp}(\Omega, N)}\delta_{\bar{s}}(f)\right\|_{\bar{p}}
. \eqno (8) %(16)
$$
for any function $f \in S_{\overline{p}, \overline{\tau}}^{\Omega}B$.

Now, since $2\le p_{j} < +\infty,$ $j=1,...,m$, using Theorem 1 from (8) we obtain
$$
\|f-S_{Q(\Omega, N)}(f)\|_{\overline{q},\overline{\theta}} <<
\left\{\sum_{\bar{s}\in \Gamma^{\perp}(\Omega, N)}\left\|\delta_{\bar{s}}(f)  \right\|_{\overline{p}}^{2}\right\}^{\frac{1}{2}} =
$$
$$
= C \left\{\sum_{\bar{s} \in \Gamma^{\perp}(\Omega, N)}\Omega^{2}(2^{-\bar{s}})
\left(\Omega^{-1}(2^{-\bar{s}})\left\|\delta_{\bar{s}}(f)  \right\|_{\overline{p}}\right)^{2}\right\}^{\frac{1}{2}}  \eqno (9)  %(17)
$$
for any function $f \in S_{\overline{p}, \overline{\tau}}^{\Omega}B$.

Item 1) proved in [4].

Let us prove item 2).
If  $\tau_{j} \le 2$, $j=1,...,m$, then using Jensen's inequality  we have
$$
\left\{\sum_{\bar{s} \in \Gamma^{\perp}(\Omega, N)}\left\|\delta_{\bar{s}}(f)  \right\|_{\overline{p}}^{2}\right\}^{\frac{1}{2}} <<
$$
$$
<< \left\|\left\{\Omega^{-1}(2^{-\bar{s}})\left\|\delta_{\bar{s}}(f)  \right\|_{\overline{p}}\right\}_{\bar{s}\in \Bbb{Z}_{+}^{m}}\right\|_{l_{\bar{\tau}}} \sup_{\bar{s}\in \Gamma^{\perp}(\Omega, N)}\Omega(2^{-\bar{s}}).
$$

Therefore, from the inequality (9) we obtain
$$
\sup \limits_{f\in S_{\overline{p}, \overline{\tau}}^{\Omega} B }
\|f-S_{Q(\Omega, N)}(f)\|_{\overline{q},\overline{\theta}} <<
\sup_{\bar{s} \in \Gamma^{\perp}(\Omega, N)}\Omega(2^{-\bar{s}}) << \frac{1}{N} ,
$$
in case $2 < p_{j} < +\infty$, $\tau_{j} \le 2$, $j=1,...,m$.
This proves the upper bound.  The lower bound in item 2) proved in [4].

Let us prove item 3).

Since $1< p_{j} \le 2,$ $j=1,...,m$, using Theorem 1 from (8) we obtain
$$
\|f-S_{Q(\Omega, N)}(f)\|_{\overline{q},\overline{\theta}} <<
\left\{\sum_{\bar{s}\in \Gamma^{\perp}(\Omega, N)}\left\|\delta_{\bar{s}}(f)  \right\|_{\overline{p}}^{p_{0}}\right\}^{\frac{1}{p_{0}}} =
$$
$$
= C \left\{\sum_{\bar{s} \in \Gamma^{\perp}(\Omega, N)}\Omega^{p_{0}}(2^{-\bar{s}})
\left(\Omega^{-1}(2^{-\bar{s}})\left\|\delta_{\bar{s}}(f)  \right\|_{\overline{p}}\right)^{p_{0}}\right\}^{\frac{1}{p_{0}}}  \eqno (10)  %(18)
$$
for any function $f \in S_{\overline{p}, \overline{\tau}}^{\Omega}B$.

If  $p_{0} < \tau_{j} < + \infty$, $j=1,...,m$, then by Holder's inequality  from (10), we get

$$
\|f-S_{Q(\Omega, N)}(f)\|_{\overline{q},\overline{\theta}} <<
\left\|\left\{\Omega^{-1}(2^{-\bar{s}})\left\|\delta_{\bar{s}}(f)  \right\|_{\overline{p}}\right\}_{\bar{s}\in \Bbb{Z}_{+}^{m}}\right\|_{l_{\bar{\tau}}}\times
$$
$$
\times \left\|\left\{\Omega(2^{-\bar{s}})\right\}_{\bar{s}\in \Gamma^{\perp}(\Omega, N)}\right\|_{l_{\bar{\epsilon}}},   \eqno (11)   %(19)
$$
where $\bar{\epsilon} = (\epsilon_{1},...,\epsilon_{m}),$ $\epsilon_{j} = 2\beta_{j}'$, $\frac{1}{\beta_{j}}+\frac{1}{\beta_{j}^{'}}=1$,
$\beta_{j} = \frac{\tau_{j}}{p_{0}}$, $j=1,...,m$.

Now by Lemma 3 and 4 from (11) we obtain
$$
\|f-S_{Q(\Omega, N)}(f)\|_{\overline{q},\overline{\theta}} << \frac{1}{N}(\log_{2}N)^{\sum\limits_{j=2}^{m}\left(\frac{1}{p_{0}}-\frac{1}{\tau_{j}}\right)}
$$
for any function $f \in S_{\overline{p}, \overline{\tau}}^{\Omega}B$.
This proves the upper bound.

Let us prove the lower bound. onsider the set similarly in \cite{ref23}
$$
\Lambda'(\Omega, N) = \Bigl\{\ \bar{s}\in \Lambda(\Omega, N) : s_{j} > \frac{1}{2ml}\log_{2}(C_{3}N), \,\, j=1,...,m\Bigr\}.
$$
N.N. Pustovoitov \cite{ref23} has been proved that, number of point is equal to $|\Lambda'(\Omega, N)| \asymp (\log_{2}N)^{m-1}$.

 After this we choose set $\bar{\Lambda(\Omega, N)}$. Lets take a number $v=[|\Lambda'(\Omega, N)|^{\frac{1}{m}}]$ - which is whole part of a number $|\Lambda'(\Omega, N)|^{\frac{1}{m}}$. Divide set $\Bbb{I}^{m} = [-\pi, \pi]^{m}$ to  $v^{m}$ ubes with side equal to $\frac{2\pi}{v}$. Then choose set $\bar{\Lambda(\Omega, N)} \subset \Lambda(\Omega, N)$, such that $|\bar{\Lambda(\Omega, N)}| = v^{m}$, , and define bijection between this set $\bar{\Lambda(\Omega, N)}$ and the set of cubes .

Let for $\bar{s}\in\bar{\Lambda(\Omega, N)}$ point $\bar{x}^{\bar{s}}$ denote centre of the cube.
Further we set notation
$$
u = \Bigl[2^{\frac{1}{m-1}\sum\limits_{j=2}^{m}(1-\frac{1}{p_{j}})(\sum\limits_{j=1}^{m}(1-\frac{1}{p_{j}}))^{-1}\log_{2}|\Lambda(\Omega, N)|} \Bigr].
$$
Consider the function %Рассмотрим функцию
$$
f_{3}(\bar{x}) = \frac{1}{N}(\log_{2}N)^{-\sum\limits_{j=2}^{m}\frac{1}{\tau_{j}}}u^{-\sum\limits_{j=1}^{m}(1-\frac{1}{p_{j}})}\Psi(\bar{x}),
$$
where  (see \cite{ref23})
$$
\Psi(\bar{x}) = \sum\limits_{\bar{s}\in \bar{\Lambda(\Omega, N)}}e^{i\langle\bar{k}^{\bar{s}}}K_{u}(\bar{x}-\bar{x}^{\bar{s}}), \;\; \bar{k}^{\bar{s}}=(k_{1}^{\bar{s}},...,k_{m}^{\bar{s}}), k_{j}^{\bar{s}} = 2^{s_{j}} + 2^{s_{j}-1}, j=1,...,m,
$$
$$
K_{u}(\bar{x}) = 2^{m}\prod_{j=1}^{m}K_{u}(x_{j}),
$$
$K_{u}(x_{j})$ -- 	is Fejer core of order $u$  by variable $x_{j}, j=1,...,m.$   Note that,
$$
u \asymp (\log_{2}N)^{\sum\limits_{j=2}^{m}(1-\frac{1}{p_{j}})(\sum\limits_{j=1}^{m}(1-\frac{1}{p_{j}}))^{-1}}. \eqno (13)
$$

In \cite{ref23} has been proved that
$$
E_{Q(\Omega, N)}(\Psi)_{1}>> |\Lambda(\Omega, N)|.  \eqno (12)
$$
Lets show that $C_{3}f_{3} \in S_{\overline{p}, \overline{\tau}}^{\Omega}B$ for some constant $C_{3} > 0.$

Since for Fejer core with one variable we have got estimation
$\|K_{u}(y)\|_{p} \asymp u^{1-\frac{1}{p}}, \,\, 1\le p \le \infty$, то
$$
\|K_{u}(\bar{x})\|_{\bar{p}} \asymp u^{\sum\limits_{j=1}^{m}(1-\frac{1}{p_{j}})}.
$$
Using this relation and $|\Lambda'(\Omega, N)| \asymp |\Lambda(\Omega, N)| \asymp (\log_{2}N)^{m-1}$ we get
$$
\left\|\left\{\Omega^{-1}(2^{-\bar{s}})\left\|\delta_{\bar{s}}(f_{3})  \right\|_{\overline{p}}\right\}_{\bar{s}\in \Bbb{Z}_{+}^{m}}\right\|_{l_{\bar{\tau}}} <<
$$
$$
<< \frac{1}{N}(\log_{2}N)^{-\sum\limits_{j=2}^{m}\frac{1}{\tau_{j}}}u^{-\sum\limits_{j=1}^{m}(1-\frac{1}{p_{j}})}\left\|\left\{\Omega^{-1}(2^{-\bar{s}})
u^{\sum\limits_{j=1}^{m}(1-\frac{1}{p_{j}})}\right\}_{\bar{s}\in \bar{\Lambda(\Omega, N)}}\right\|_{l_{\bar{\tau}}} <<
$$
$$
<< (\log_{2}N)^{-\sum\limits_{j=2}^{m}\frac{1}{\tau_{j}}}\left\|\left\{1\right\}_{\bar{s}\in \bar{\Lambda(\Omega, N)}}\right\|_{l_{\bar{\tau}}}.
$$
Since by Lemma 4  estimation
$$
\left\|\left\{1\right\}_{\bar{s}\in \bar{\Lambda(\Omega, N)}}\right\|_{l_{\bar{\tau}}}<<(\log_{2}N)^{\sum\limits_{j=2}^{m}\frac{1}{\tau_{j}}}
$$
is true, then
$$
\left\|\left\{\Omega^{-1}(2^{-\bar{s}})\left\|\delta_{\bar{s}}(f)  \right\|_{\overline{p}}\right\}_{\bar{s}\in \Bbb{Z}_{+}^{m}}\right\|_{l_{\bar{\tau}}} << C.
$$
Because $C_{3}f_{3} \in S_{\overline{p}, \overline{\tau}}^{\Omega}B$.

Since $L_{\bar{q}, \bar{\theta}}(\Bbb{I}^{m}) \subset L_{1}(\Bbb{I}^{m})$ and $\|f\|<< \|f\|_{\bar{q}, \bar{\theta}},$ then
$$
E_{Q(\Omega, N)}(f_{3})_{_{\bar{q}, \bar{\theta}}} >> E_{Q(\Omega, N)}(f_{3})_{1} = C\frac{1}{N}(\log_{2}N)^{-\sum\limits_{j=2}^{m}\frac{1}{\tau_{j}}}u^{-\sum\limits_{j=1}^{m}(1-\frac{1}{p_{j}})}E_{Q(\Omega, N)}(\Psi)_{1}.
$$
Therefore, by the estimates (12), (13) we have got
$$
E_{Q(\Omega, N)}(f_{3})_{_{\bar{q}, \bar{\theta}}} >> \frac{1}{N}(\log_{2}N)^{-\sum\limits_{j=2}^{m}\frac{1}{\tau_{j}}}u^{-\sum\limits_{j=1}^{m}(1-\frac{1}{p_{j}})}|\Lambda(\Omega, N)| >>
$$
$$
>> \frac{1}{N}(\log_{2}N)^{-\sum\limits_{j=2}^{m}\frac{1}{\tau_{j}}}(\log_{2}N)^{-\sum\limits_{j=2}^{m}(1-\frac{1}{p_{j}})}(\log_{2}N)^{m-1} = C\frac{1}{N}(\log_{2}N)^{\sum\limits_{j=2}^{m}(\frac{1}{p_{j}}-\frac{1}{\tau_{j}})}.
$$
The Theorem 5 is proved.

Now consider the case $q_{j} = p_{j}, \;\; j=1,...,m$ and  $\Omega(\bar{t}) = \prod\limits_{j=1}^{m}t_{j}^{r_{j}}, \,\, r_{j} >0, \,\, t_{j} \in [0, 1], j=1,...m.$

\begin{theorem} %Теорема 6
Let
$\bar{r} = (r_{1},...,r_{m}), \,\, 0 < r_{1}=...=r_{\nu}< r_{\nu+1}\le...\le r_{m}$ and $1 < q_{j} < p_{j} < \infty$, $p_{j} \ge 2$, $1< \theta_{j}< \infty$, $1\le \tau_{j} \le +\infty,$ $j=1,...,m$.

 If $2\le p_{j} < \theta_{j} < \infty$, $2 \le \tau_{j} < +\infty,$
$j=1,...,m,$  then
$$
\sup \limits_{f\in
S_{\overline{p}, \overline{\theta}, \overline{\tau}
}^{\bar{r}}B }\|f-S_{N}^{\gamma}(f)\|_{\overline{p}} << N^{-r_{1}}\left(\log_{2}N\right)^{\sum\limits_{j=1}^{m}
\left(\frac{1}{p_{j}}-\frac{1}{\theta_{j}}\right)}  \left(\log_{2}N\right)^{\sum\limits_{j=2}^{m}
\left(\frac{1}{2}-\frac{1}{\tau_{j}}\right)}
$$
and if $p_{1} = ... =p_{m}=p,$ then
$$
N^{-r_{1}}\left(\log_{2}N\right)^{\sum\limits_{j=1}^{m}
\left(\frac{1}{p}-\frac{1}{\theta_{j}}\right)}  \left(\log_{2}N\right)^{\sum\limits_{j=2}^{m}
\left(\frac{1}{p}-\frac{1}{\tau_{j}}\right)} << \sup \limits_{f\in
S_{p, \overline{\theta}, \overline{\tau}
}^{\bar{r}}B }\|f-S_{N}^{\gamma}(f)\|_{\overline{p}}
$$
\end{theorem}

{\bf Proof.} Let $f \in S_{\bar{p}, \overline{\theta}, \overline{\tau}}^{\bar{r}}B$.
Now, since $2\le p_{j} < +\infty,$ $2 \le \tau_{j} < +\infty,$ $j=1,...,m$, using Theorem 1 and the inequality of different metric for trigonometric polynomials (see \cite{ref05}), the inequality Holder's  we obtain
$$
\|f\|_{\overline{p}} <<
\left\{\sum_{\bar{s}\in \Bbb{Z}^{m}}\left\|\delta_{\bar{s}}(f)\right\|_{\overline{p}}^{2}\right\}^{\frac{1}{2}} <<
$$
$$
<< \left\{\sum_{\bar{s}\in \Bbb{Z}^{m}}\left\|\delta_{\bar{s}}(f)\right\|_{\overline{p}, \bar\theta}^{2}\prod\limits_{j=1}^{m}(s_{j}+1)^{\frac{1}{p_{j}}-\frac{1}{\theta_{j}}}\right\}^{\frac{1}{2}} <<
$$
$$
<< \left\|\left\{2^{\langle\bar{s}, \bar{r}\rangle}\left\|\delta_{\bar{s}}(f)\right\|_{\overline{p}}\right\}_{\bar{s}\in \Bbb{Z}_{+}^{m}}\right\|_{l_{\bar{\tau}}}
 \left\|\left\{2^{-\langle\bar{s}, \bar{r}\rangle}\prod\limits_{j=1}^{m}(s_{j}+1)^{\frac{1}{p_{j}}-\frac{1}{\theta_{j}}}\right\}_{\bar{s}\in \Bbb{Z}_{+}^{m}}\right\|_{l_{\bar{\epsilon}}},   \eqno (14) %(12)   %(20)
$$
where $\bar{\epsilon} = (\epsilon_{1},...,\epsilon_{m}),$ $\epsilon_{j} = \frac{2\tau_{j}}{\tau_{j}-2}$,  $j=1,...,m$.
Taking into account that $r_{j} > 0, \,\, j=1,...,m$  we get
$$
\left\|\left\{2^{-\langle\bar{s}, \bar{r}\rangle}\prod\limits_{j=1}^{m}(s_{j}+1)^{\frac{1}{p_{j}}-\frac{1}{\theta_{j}}}\right\}_{\bar{s}\in \Bbb{Z}_{+}^{m}}\right\|_{l_{\bar{\epsilon}}} < \infty.
$$
Hence, it follows from (14) that $S_{\bar{p}, \overline{\theta}, \overline{\tau}}^{\bar{r}}B \subset L_{\bar{p}}(\Bbb{I}^{m})$ and
$$
\|f - S_{N}^{\gamma}(f)\|_{\overline{p}} <<
\left\|\left\{2^{\langle\bar{s}, \bar{r}\rangle}\left\|\delta_{\bar{s}}(f)  \right\|_{\overline{p}}\right\}_{\bar{s}\in \Bbb{Z}_{+}^{m}}\right\|_{l_{\bar{\tau}}}
 \left\|\left\{2^{-\langle\bar{s}, \bar{r}\rangle}\prod\limits_{j=1}^{m}(s_{j}+1)^{\frac{1}{p_{j}}-\frac{1}{\theta_{j}}}\right\}_{\bar{s}\in \Gamma^{\perp}(N)}\right\|_{l_{\bar{\epsilon}}}   \eqno (15)  %(13)   %(21)
 $$
where $\Gamma^{\perp}(N) = \{\bar{s}\in \Bbb{Z}_{+}^{m}: \,\, \langle\bar{s}, \bar{\gamma}\rangle \ge \log_{2}N^{\frac{1}{r_{1}}}\}$.

Next applying inequality
$$
I_{N} = \left\|\left\{2^{-\langle\bar{s}, \bar{\gamma}\rangle\beta}\prod\limits_{j=1}^{m}s_{j}^{d_{j}} \right\}_{\bar{s} \in \Gamma^{\perp}(N)}\right\|_{l_{\bar{\theta}}} <<2^{-n\beta}n^{\sum_{j=1}^{m}d_{j} + \sum_{j=2}^{m}\frac{1}{\theta_{j}}}
$$
for $\beta > 0,$ $d_{j} \ge 0, j=1,...,m$, then
$$
\left\|\left\{2^{-\langle\bar{s}, \bar{\gamma}\rangle r_{1}}\prod\limits_{j=1}^{m}(s_{j}+1)^{\frac{1}{p_{j}}-\frac{1}{\theta_{j}}}\right\}_{\bar{s}\in \Gamma^{\perp}(N)}\right\|_{l_{\bar{\epsilon}}}<<
$$
$$
<< N^{-r_{1}}(\log_{2}N)^{\sum_{j=1}^{m}(\frac{1}{p_{j}}-\frac{1}{\theta_{j}}) + \sum_{j=2}^{m}\frac{1}{\varepsilon_{j}}}.
$$
Therefore from (15) we obtain
$$
\|f - S_{N}^{\gamma}(f)\|_{\overline{p}} << N^{-r_{1}}(\log_{2}N)^{\sum_{j=1}^{m}(\frac{1}{p_{j}}-\frac{1}{\theta_{j}}) + \sum_{j=2}^{m}(\frac{1}{2} - \frac{1}{\tau_{j}})}
$$
for any function $f \in S_{\overline{p}, \bar\theta, \overline{\tau}}^{\bar{r}}B$.
This proves the upper bound.

Let us prove the lower bound. Consider the function
$$
f_{4}(\bar{x})= (\log_{2}N)^{-\sum\limits_{j=2}^{m}\frac{1}{\tau_{j}}}
\sum_{\langle\bar{s}, \bar{\gamma}\rangle = \log_{2}N^{\frac{1}{r_{1}}}}\prod_{j=1}^{m}2^{-s_{j}r_{j}}s_{j}^{-\frac{1}{\theta_{j}}}\sum_{\bar{k}\in\rho(\bar
s)}\prod_{j=1}^{m}(k_{j} - 2^{s_{j}-1}+1)^{\frac{1}{p}-1} e^{i\langle\bar{k},\bar{x}\rangle}.
$$
Then $f_{4} \in L_{\bar{p}, \bar{\theta}}(\Bbb{I}^{m})$. Now , by relation
$$
\Bigl\|\sum_{\bar{k}\in\rho(\bar
s)}\prod_{j=1}^{m}(k_{j} - 2^{s_{j}-1}+1)^{\frac{1}{p}-1} e^{i\langle\bar{k},\bar{x}\rangle}\Bigr\|_{\bar{p}, \bar{\theta}}\asymp
\prod_{j=1}^{m}(s_{j}+1)^{\frac{1}{\theta_{j}}}  \eqno (16)  %(14)  %(23)
$$
for $1 < p_{j} < \infty, \,\, 1<\theta_{j}< \infty, j=1,...,m$ and by Lemma 1 \cite{ref04} we get
$$
\left\|\left\{2^{\langle\bar{s}, \bar{r}\rangle}\left\|\delta_{\bar{s}}(f_{4})\right\|_{\bar{p}, \bar{\theta}}\right\}_{\bar{s}\in \Bbb{Z}_{+}^{m}}\right\|_{l_{\bar{\tau}}} << (\log_{2}N)^{-\sum\limits_{j=2}^{m}\frac{1}{\tau_{j}}}\left\|\left\{1\right\}_{\bar{s}\in \varkappa(N)}\right\|_{l_{\bar{\tau}}} \le C,
$$
where $\varkappa(N) = \{\bar{s}\in \Bbb{Z}_{+}^{m} : \,\, \langle\bar{s}, \bar{\gamma}\rangle = \log_{2}N^{\frac{1}{r_{1}}} \}$.

Hence the function $C_{4}f_{4} \in S_{\bar{p}, \bar\theta, \overline{\tau}}^{\bar{r}}B$.

Since $2\le p = p_{1}=...= p_{m} < \infty$, then by Littlewood-Paley theorem \cite{ref21} we obtain
$$
\|f_{4}-S_{N}^{\gamma}(f_{4})\|_{p} = \|f_{4}\|_{p} >> \Bigl\|\Bigl\{\sum_{\bar{s} \in \varkappa(N)} |\delta_{\bar{s}}(f_{4})|^{2}\Bigr\}^{\frac{1}{2}}   \Bigr\|_{p} >> \Bigl\{\sum_{\bar{s} \in \varkappa(N)}\left\|\delta_{\bar{s}}(f_{2})\right\|_{p}^{p}\Bigr\}^{\frac{1}{p}}.
$$
By relation (16) for $\theta_{j} = p_{j} = p, j=1,...,m$ , it follows that
$$
\|f_{4}-S_{N}^{\gamma}(f_{4})\|_{p} >> (\log_{2}N)^{-\sum\limits_{j=2}^{m}\frac{1}{\tau_{j}}}\Bigl(\sum_{\bar{s} \in \varkappa(N)}2^{-\langle\bar{s}, \bar{r}\rangle p}\prod_{j=1}^{m}(s_{j}+1)^{(\frac{1}{p} - \frac{1}{\theta_{j}})p} \Bigr)^{\frac{1}{p}} >>
$$
$$
>> N^{-r_{1}}(\log_{2}N)^{\sum\limits_{j=1}^{m}(\frac{1}{p} -\frac{1}{\theta_{j}})}(\log_{2}N)^{\sum\limits_{j=2}^{m}(\frac{1}{p} -\frac{1}{\tau_{j}})}
$$
for  function $C_{4}f_{4} \in S_{p, \bar\theta, \overline{\tau}}^{\bar{r}}B$. So Theorem 6 has been proved.

{\bf Remark.} Note that for the case  $q_{j}=\theta_{j}=q$, $p_{j}=p$, $\tau_{j} = \tau$, $j=1,...,m$, Theorem 5 was proved by S.A. Stasyuk \cite{ref30}. For the case  $p_{j}=\theta_{j}^{(1)}=p$, $q_{j}= \theta_{j}^{(2)} = q$, $\tau_{j} = +\infty$, $j=1,...,m$, Theorem 4 was proved by N.N. Pustovoitov \cite{ref23}.

{\bf Acknowledgements.}  This work was supported by the Ministry of Education and Science of Republic Kazakhstan (Grant no. 5129GF4) and by the Competitiveness Enhancement Program of the Ural Federal University (Enactment  of the Government of the Russian Federation of March 16, 2013 no. 211, agreement no. 02.A03. 21.0006 of August 27, 2013).

\end{document}